\documentclass[12pt]{article}

\usepackage{amssymb,amsthm,amscd,array}

\let\ssection=\section
\renewcommand{\section}{\setcounter{equation}{0}\ssection}

\newcommand{\bbR}{\mathbb{R}}

\newcommand{\bbC}{\mathbb{C}}

\newcommand{\bbN}{\mathbb{N}}

\newcommand{\cC}{{\mathcal{C}}}
\newcommand{\rC}{{\mathrm{C}}}

\newcommand{\cD}{{\mathcal{D}}}
\newcommand{\ce}{\mathrm{ce}}
\newcommand{\cE}{{\mathcal{E}}}
\newcommand{\diff}{\mathrm{diff}}
\newcommand{\Diff}{\mathrm{Diff}}
\newcommand{\Div}{\mathrm{Div}}
\newcommand{\rD}{\mathrm{D}}
\newcommand{\rd}{\mathrm{d}}
\newcommand{\rE}{\mathrm{E}}
\newcommand{\Euler}{{\mathcal E}}
\newcommand{\End}{\mathrm{End}}
\newcommand{\cF}{{\mathcal{F}}}

\newcommand{\rg}{\mathrm{g}}

\newcommand{\rG}{\mathrm{G}}

\newcommand{\rh}{\mathrm{h}}
\newcommand{\cI}{{\mathcal{I}}}
\newcommand{\Id}{\mathrm{Id}}
\newcommand{\cL}{\mathcal{L}}

\newcommand{\cN}{\mathcal{N}}
\newcommand{\rN}{\mathrm{N}}
\newcommand{\Pol}{\mathrm{Pol}}
\newcommand{\cQ}{{\mathcal{Q}}}
\newcommand{\rR}{\mathrm{R}}

\newcommand{\cS}{{\mathcal{S}}}

\newcommand{\SL}{\mathrm{SL}}
\newcommand{\Sl}{\mathrm{sl}}
\newcommand{\SO}{\mathrm{SO}}
\newcommand{\se}{\mathrm{e}}
\newcommand{\so}{\mathrm{o}}
\newcommand{\Sp}{\mathrm{Sp}}

\newcommand{\Symm}{\mathrm{Symm}}
\newcommand{\rT}{\mathrm{T}}
\newcommand{\Tr}{\mathrm{Tr}}

\newcommand{\Vect}{\mathrm{Vect}}

\newcommand{\Vol}{\mathrm{Vol}}
\newcommand{\rW}{\mathrm{W}}
\newcommand{\cZ}{{\mathcal{Z}}}

\newcommand{\half}{\frac{1}{2}}
\newcommand{\fg}{\mathfrak{g}}

\def\mod#1{\left|{#1}\right|}

\begin{document}



\def\a{\alpha}
\def\b{\beta}
\def\d{\delta}
\def\g{\gamma}
\def\om{\omega}
\def\r{\rho}
\def\s{\sigma}
\def\vfi{\varphi}
\def\vr{\varrho}
\def\l{\lambda}
\def\m{\mu}
\def\implies{\Rightarrow}

\oddsidemargin .1truein
\newtheorem{thm}{Theorem}[section]
\newtheorem{lem}[thm]{Lemma}
\newtheorem{cor}[thm]{Corollary}
\newtheorem{pro}[thm]{Proposition}
\newtheorem{ex}[thm]{Example}
\newtheorem{rmk}[thm]{Remark}
\newtheorem{defi}[thm]{Definition}


\title{Conformally equivariant quantization:\\
Existence and uniqueness}

\author{C.~Duval\thanks{
Universit\'e de la M\'editerran\'ee and CPT-CNRS,
Luminy Case 907, F--13288 Marseille, Cedex~9, FRANCE;
mailto:duval@cpt.univ-mrs.fr
}
\and
P.~Lecomte\thanks{
Institut de Math\'ematiques, Universit\'e de Li\`ege, Sart Tilman, Gde
Traverse, 12 (B~37), B--4000 Li\`ege, BELGIUM;
mailto:plecomte@ulg.ac.be
}
\and
V.~Ovsienko\thanks{
CNRS, Centre de Physique Th\'eorique,
CPT-CNRS, Luminy Case 907,
F--13288 Marseille, Cedex~9, FRANCE;
mailto:ovsienko@cpt.univ-mrs.fr
}
}

\date{}

\maketitle

\hfill{\textit{\small In Memory of
Moshe Flato 
and
Andr\'e Lichnerowicz}}

\bigskip

\thispagestyle{empty}

\begin{abstract}
We prove the existence and the uniqueness of a conformally equivariant symbol calculus
and quantization on any conformally flat pseudo-Riemannian manifold~$(M,\rg)$. In other
words, we establish a canonical isomorphism between the spaces of polynomials on~$T^*M$
and of differential operators on tensor densities over~$M$, both viewed as modules over
the Lie algebra $\so(p+1,q+1)$ where $p+q=\dim(M)$. This quantization exists for generic
values of the weights of the tensor densities and compute the critical values of the
weights yielding obstructions to the existence of such an isomorphism.
In the particular case of half-densities, we obtain a conformally invariant
star-product.
\end{abstract}


\newpage

\section{Introduction}

The general problem of quantization is often understood as the quest for a
correspondence between the smooth functions of a given symplectic manifold, i.e., the
classical observables, and (symmetric) operators on a certain associated Hilbert space,
which are called the quantum observables. This correspondence must satisfy a number of
additional properties that heavily depend upon the standpoint of the authors and are,
by no means, universal.

One guiding principle for the search of a quantization procedure is to impose
further coherence with some natural symmetry of phase space. This
constitutes the foundations of the ``orbit method'' \cite{Kir1},  geometric
quantization \cite{Kos1,JMS,Kir3} in the presence of symmetries, Moyal-Weyl
quantization (see, e.g.,~\cite{Fed}) defined by requiring invariance with respect to
the linear symplectic group $\Sp(2n,\bbR)$ of~$\bbR^{2n}$.

\subsection{Equivariant quantization problem}

In all previous examples, the symmetry group was a Lie subgroup of the group of all
symplectomorphisms of the symplectic manifold. Here, we will confine considerations to
the case of cotangent bundles, $T^*M$, with their canonical (vertical) polarization.
This polarization, together with the Liouville $1$-form on $T^*M$, should be preserved
by the symmetry group which naturally arises as the cotangent lift of a Lie group, $G$,
acting on $M$.

We will thus look for an identification, as $G$-modules, between the space, $\cS(M)$, of
smooth functions on $T^*M$ that are polynomial on the fibers and the the space,
$\cD(M)$, of linear differential operators on~$M$.

Note that the Moyal-Weyl quantization does not fit into this general framework in which
the symmetries of configuration space, $M$, play a central r\^ole. Indeed, the action of
$\Sp(2n,\bbR)$ on $T^*\bbR^n$ does not descend to $\bbR^n$.

Now, $G$-equivariance between $\cS(M)$ and $\cD(M)$ is clearly too strong a
requirement if $G$ is the group, $\Diff(M)$, of all diffeomorphisms of $M$. We will
therefore impose such an equivariance in the weaker case where $G$ is a
finite-dimensional Lie group whose action on $M$ is to be only local.

The main tool we will be using is provided by the notion of (flat) $G$-structure.
Let us recall that a $G$-structure on $M$ is defined by a local action of $G$ on
$M$, compatible with a local identification of~$M$ with some homogeneous space~$G/H$. 
More precisely, it is defined by an atlas of charts $(\varphi_\a,V_\a)$ with
$\varphi_\a:V_\a\to{}G/H$ such that $\varphi_\a\circ\varphi_\b^{-1}$ is given by an
element of $G$.
This new approach significantly differs from the more usual one which makes use of
connections to intrinsically define quantization procedures and symbol calculus.

The $G$-structure we will consider in this article is the conformal structure
with $G=\SO(p+1,q+1)$ modeled on the pseudo-Riemannian manifold $S^p\times{}S^q$.

Our purpose is to show that there exists, actually, a canonical isomorphism of
$\SO(p+1,q+1)$-modules between the space of symbols, $\cS(M)$, and the space of
differential operators, $\cD(M)$. 

Experience of other approaches to the quantization problem and of the geometrical study
of differential equations prompts us to rather consider the space, $\cD_{\l,\m}(M)$,
of differential operators with arguments and values in the space of tensor densities of
weights $\l$ and $\m$ respectively.  
So, we will naturally need to study this space of differential operators as a
$\SO(p+1,q+1)$-module. As a consequence, the $\SO(p+1,q+1)$-module of symbols will be
twisted by the weight $\d=\m-\l$, and denoted by $\cS_\d(M)$.

Let us emphasize that equivariant quantization of $G$-struc\-tures has
already been carried out in the case of projective structures,
i.e.~$\SL(n+1,\bbR)$-structures, in the recent papers \cite{LO,Lec2}. As for conformal
structures, a first step towards their equivariant quantization was taken in
\cite{DO2} in the case of second order operators.

In the particular case $n=1$, both conformal and projective structures coincide. We
refer to \cite{CMZ} for a thorough study of $\Sl(2,\bbR)$-equivariant
quantization and of the corresponding invariant star-product. See also
\cite{Wil} for a classic monography on the structures of $\Sl(2,\bbR)$-module on the
space of differential operators on the real line.

There exist various approaches to the quantization problem, however, our viewpoint put
emphasis on the equivariance condition with respect to a (maximal) group, $G$, of
symmetry in the context of deformation quantization. Now, $G$-equivariance is the
root of geometric quantization \cite{JMS,Kos1,Kir3}, Berezin quantization
\cite{Ber1,Ber2}, etc., but it seems to constitute a fairly new approach in the
framework of symbol calculus, deformation theory and semi-classical approximations
dealt with in this work.

\subsection{Quantizing equivariantly conformal structures}

We outline here the main results we have obtained, and describe the general
framework adopted in this article to answer the question raised in the
preceding section.

We review in Section \ref{BDT} the structures of the spaces of symbols $\cS_\d(M)$ and
of differential operators $\cD_{\l,\m}(M)$ as $\Vect(M)$-modules. The Lie algebra,
$\so(p+1,q+1)$, of conformal Killing vectors of a conformally flat manifold $(M,\rg)$
is then described. The restriction of the preceding $\Vect(M)$-actions is also
explicitly calculated.

Section \ref{MainResults} presents the main theorems which establish the existence and
the uniqueness of a $\so(p+1,q+1)$-equivariant quantization in the special and
fundamental case $\l=\m$. It turns out that the value $\l=\m=\half$ guarantees that our
quantization actually defines a star-product on $T^*M$. This is precisely the value of
the weights used in geometric quantization. 
In the general case, we again obtain a canonical isomorphism of $\so(p+1,q+1)$-modules,
except for an infinite series of values of $\d=\m-\l$, which we call resonances. 

\goodbreak

Section \ref{ConfInvOp} is devoted to 
the algebra of invariants.
One considers the action of $\so(p+1,q+1)$ and the Euclidean subalgebra $\se(p,q)$ on
the space of polynomials on $T^*\bbR^n$. Resorting to Weyl's theory of
invariants, we characterize the commutant of $\so(p+1,q+1)$ within the algebra of
operators on the latter space of polynomials: it is a commutative associative algebra
with two generators. The commutant of $\se(p,q)$ has already been
determined in
\cite{DO2}. These algebras of conformal and Euclidean invariants play a crucial r\^ole
in our work and enable us to compute the Casimir operators $\rC_\d$ and $\cC_{\l,\m}$
of the $\so(p+1,q+1)$-actions on $\cS_\d(M)$ and $\cD_{\l,\m}(M)$ respectively.

In Section \ref{ProofMain} we provide the proofs of the main theorems.
It should be stressed that these proofs rely on the diagonalization of the
preceding Casimir operators in an essential way. The same idea has already been
exploited in \cite{CMZ} in the one-dimensional case, and in \cite{Lec2} in the case of
projectively flat manifolds of dimension~$n>1$.

Section \ref{quadratic} is concerned with the explicit expression of the quantization
map restricted to second order polynomials on $T^*M$. It is worth noticing that our
conformally equivariant quantization on $T^*\bbR^n$ differs from the standard Weyl
quantization to which it constitutes a new alternative.

This article is also related to various different subjects, namely representation
theory, the theory of invariant differential operators and the cohomology of Lie
algebras. A number of very concrete problems could be tackled in this framework. For
example, the quantization of the geodesic flow has been achieved in a purely
conformally invariant manner \cite{DO2}. Also the Yamabe-Laplace operator (or conformal
Laplacian), see \cite{Bes}, arose from the quantization of the same geodesic flow in a
resonant case (recall that this operator is of special importance in field
theory in a curved space-time, see, e.g., \cite{PR}).

Let us finally mention that this work opens up a number of original questions under
current investigation, viz the determination of the $\so(p+1,q+1)$-invariant
star-product and multi-dimensional Schwarzian derivative.

\medskip
\textbf{Acknowledgments:} We are indebted to A.~A.~Kirillov for most enlightening
discussions and also to S.~Loubon-Djounga for his efficient help. Special thanks are
due to D.~Leites for clarifying conversations.

\newpage

\section{Basic definitions and tools}\label{BDT}

\subsection{Differential operators on tensor densities}

Let us recall that a tensor density of degree~$\l$ on manifold $M$ is a smooth
section of the line bundle $\Delta_\l(M)=\mod{\Lambda^nT^*M}^{\otimes\l}$ over $M$.
The space of tensor densities of degree~$\l$ is naturally a $\Diff(M)$- and
$\Vect(M)$-module. In this paper, we will consider the space $\cF_\l(M)$ (or $\cF_\l$ in
short) of complex-valued smooth tensor densities, i.e., of the sections of
$\Delta_\l(M)\otimes\bbC$.

The space $\cD_{\l,\m}$ of linear differential operators
\begin{equation}
A:\cF_\l\to\cF_\m
\label{A}
\end{equation}
from $\l$-densities to $\m$-densities on $M$ is naturally a
$\Diff(M)$- and $\Vect(M)$-module.
These modules have been studied and classified in
\cite{DO1,LMT,LO,GO,Gar,Mat,DO2,Lec1}.

There is a filtration
$\cD^0_{\l,\m}\subset\cD^1_{\l,\m}\subset\cdots\subset\cD^k_{\l,\m}\subset\cdots$,
where the module of zero-order operators $\cD^0_{\l,\m}\cong\cF_{\m-\l}$
consists of multiplication by $(\m-\l)$-densities. The higher-order
modules are
defined by induction: $A\in\cD^k_{\l,\m}$ if $[A,f]\in\cD^{k-1}_{\l,\m}$
for
every $f\in{}C^\infty(M)$.

\subsection{Symbols with values in tensor densities}\label{SymbDens}

Consider the space $\cS=\Gamma(S(TM))$ of contravariant symmetric tensor
fields on~$M$ which is naturally a $\Diff(M)$- and $\Vect(M)$-module. We will define
the space of $\d$-weighted symbols on $T^*M$ as the space of sections
\begin{equation}
\cS_\d
=
\Gamma\left(
S(TM)\otimes\Delta_\d(M)
\right).
\label{P}
\end{equation}
The space $\cS_\d$ is also, naturally, a $\Diff(M)$- and $\Vect(M)$-module.

Again, there is a filtration
$\cS^0_\d\subset\cS^1_\d\subset\cdots\subset\cS^k_\d\subset\cdots$,
where $\cS^k_\d$ denotes the space of symbols of degree less or equal to $k$. 
In contrast to the filtration on the space $\cD_{\l,\m}$ of differential operators, the
above filtration on the space (\ref{P}) of symbols actually leads to a
$\Diff(M)$-invariant graduation 
\begin{equation}
\cS_\d=\bigoplus_{k=0}^\infty{\cS_{k,\d}}
\label{grad}
\end{equation}
where $\cS_{k,\d}$ denotes the space of homogeneous polynomials
(isomorphic to $\cS^k_\d/\cS^{k-1}_\d$).
\goodbreak

\subsection{$\cD_{\l,\m}$ and $\cS_\d$ as $\Vect(M)$-modules}\label{VectModules}

We will always assume $M$ orientable and identify $\cF_\l$ to $C^\infty(M)\otimes\bbC$
by the choice of a volume form, $\Vol$, on $M$.

It is clear from the definition of $\cD_{\l,\m}$ that the
corresponding $\Vect(M)$-action, $\cL^{\l,\m}$, is given by
\begin{equation}
\cL^{\l,\m}_X(A) = 
L^{\m}_X\,A - A\,L^{\l}_X
\label{Dinf}
\end{equation}
where $X\in\Vect(M)$, and $L^\l_X$ is the standard Lie derivative of $\l$-densities
$\cF_\l$. Now, any $\l$-density being represented
by $f\,|\Vol|^\l$ for some $f\in{}C^\infty(M)\otimes\bbC$, the Lie derivative~$L^\l_X$
is thus given by
\begin{equation}
L^\l_X(f) = 
X(f)+\l\,\Div(X)\,f
\label{Finf}
\end{equation}
where $\Div(X)=L_X(\Vol)/\Vol$.

As to the $\Vect(M)$-action, $L^\d$, on $\cS_\d$, it reads
\begin{equation}
L_X^\d(P)=
L_X(P)+\d\,\Div(X)\,P
\label{Sinf}
\end{equation}
where $L_X$ denotes here the Lie derivative of contravariant tensors given by the
cotangent lift of $X\in\Vect(M)$.

\subsection{The modules $\cD_{\l,\m}(\bbR^n)$ and $\cS_\d(\bbR^n)$}\label{TheModules}

In a given coordinate system~$(x^1,\ldots,x^n)$ on $\bbR^n$, the expression of a
differential operator $A\in\cD^k_{\l,\m}$ (see (\ref{A})) reads
\begin{equation}
A
=
A_k^{{i_1}\ldots{i_k}}\partial_{i_1}\ldots\partial_{i_k}
+
\cdots
+
A_1^i\partial_i
+
A_0
\label{DiffOp}
\end{equation}
where $\partial_i=\partial/\partial x^i$, the coefficient
$A_\ell^{{i_1}\ldots{i_\ell}}\in{}C^\infty(\bbR^n)$ being symmetric in
${i_1},\dots,{i_\ell}$ for
$\ell=0,1,\ldots,k$.
{}From now on we suppose a summation over repeated indices.

The local expression of a symbol $P\in\cS^k_\d$ (see (\ref{P})), in the canonical
coordinate system $(x^1,\ldots,x^n,\xi_1,\ldots,\xi_n)$ on $T^*\bbR^n$ is then
\begin{equation}
P
=
P_k^{{i_1}\ldots{i_k}}\xi_{i_1}\ldots\xi_{i_k}
+
\cdots
+
P_1^i\xi_i
+
P_0
\label{Symbol}
\end{equation}
where $P_\ell^{{i_1}\ldots{i_\ell}}\in{}C^\infty(\bbR^n)$ represent the
components of symmetric contravariant tensor fields on $\bbR^n$ (for
$\ell=0,1,\ldots,k$).

As vector spaces, $\cD_{\l,\m}$ and $\cS_\d$ are clearly isomorphic, though not in a
canonical way. 
For example, the normal ordering map
\begin{equation}
\s:A_k^{{i_1}\ldots{i_k}}
\partial_{i_1}\cdots
\partial_{i_k}
\mapsto
A_k^{{i_1}\ldots{i_k}}
\xi_{i_1}\cdots
\xi_{i_k}
\label{sigma}
\end{equation}
defines such an isomorphism.

\goodbreak

The $\Vect(\bbR^n)$-action (\ref{Dinf}) on differential operators is, of course,
different from the standard $\Vect(\bbR^n)$-action (\ref{Sinf}) on polynomials.
We will, therefore, distinguish the two $\Vect(\bbR^n)$-modules
\begin{eqnarray}
\label{DvaMod}
\cD_{\l,\m}
&=&
(\Pol(T^*\bbR^n),\cL^{\l,\m}),\\
\cS_\d
&=&
(\Pol(T^*\bbR^n),L^\d).
\end{eqnarray}
In particular, a vector field $X$ corresponds to a first-order polynomial,
$X=X^i\xi_i$. The operator of Lie derivative is then given by the Hamiltonian
vector field
\begin{equation}
L^\d_X=
\partial_{\xi_i}X\partial_i
-
\partial_iX\partial_{\xi_i}
+\d\,\rD{X},
\label{LieDer1}
\end{equation}
where $\rD=\partial_{\xi_i}\partial_i$ is the divergence operator (see Section
\ref{Euclinv}). This local expression precisely corresponds to the previous
expression (\ref{Sinf}).

One easily proves the
\begin{pro}
The $\Vect(\bbR^n)$-action on $\cD_{\l,\m}$ has the following form
\begin{equation}
\begin{array}{ccl}
\cL_X^{\l,\m}&=&
L^\d_X-\frac{1}{2}\,\partial_i\partial_jX
\partial_{\xi_i}\partial_{\xi_j}-
\l\,(\partial_i\circ\rD)X
\partial_{\xi_i}\\[8pt]
&&+(\hbox{higher order derivatives
$\partial_{i_1}\cdots\partial_{i_\ell}X$})
\end{array}
\label{LieDer2}
\end{equation}
where $\d=\m-\l$.
\label{PolPro}
\end{pro}

\subsection{Conformally flat manifolds}

Throughout this paper we will deal with conformally flat manifolds. Let us recall
that a smooth pseudo-Riemannian manifold $(M,\rg)$ is conformally flat if,
for every $x\in{}M$, there exists a neighborhood $V_x$ of $x$ and 
$F\in{}C^\infty(V_x,\bbR^*_+)$ such that $(V_x,g)$ is flat with the new metric
$g=F\,\rg$.

The basic example of a Riemannian $n$-dimensional conformally flat manifold is the
sphere $S^n$ with its canonical metric, and $S^p\times S^q$ in the case of
signature~$p-q$.  A conformally flat manifold is locally identified with such a
homogeneous space and thus admits a local action of $\SO(p+1,q+1)$. The associated,
locally defined, action of the Lie algebra $\so(p+1,q+1)$ corresponds (if
$n=p+q\geq3$) to that of the subalgebra of the vector fields $X$ solutions of
$$
L_X\rg=f\,\rg
$$
for some $f\in{}C^\infty(M)$ depending upon $X$. 


It is well known that a conformally flat manifold
admits an atlas in which $\so(p+1,q+1)$ is generated by
\begin{equation}
\matrix{
X_i  
&=&
\displaystyle \frac{\partial}{\partial x^i}\;\hfill\cr
\noalign{\smallskip}
X_{ij}
&=&
\displaystyle x_i\frac{\partial}{\partial x^j}-
x_j\frac{\partial}{\partial x^i}\;\hfill\cr
\noalign{\smallskip}
X_0  
&=&
\displaystyle x^i\frac{\partial}{\partial x^i}\;\hfill\cr
\noalign{\smallskip}
\bar X_i 
&=&
\displaystyle x_jx^j\frac{\partial}{\partial x^i}-
2x_ix^j\frac{\partial}{\partial x^j}\hfill\cr
}
\label{confGenerators}
\end{equation}
where $i,j=1,\ldots,n$ and $x_i=g_{ij}x^j$. 
In the sequel, indices will be
raised and lowered by means of the (flat) metric~$g$.

Let us introduce the following nested Lie subalgebras that will be considered below,
namely
\begin{equation}
\so(p,q)\subset\se(p,q)\subset\ce(p,q)\subset\so(p+1,q+1)
\label{subalgebras}
\end{equation}
where $\so(p,q)$ is generated by the $X_{ij}$, the Euclidean
sub\-algebra $\se(p,q)$ by $X_{ij}$ and $X_i$ and the Lie algebra
$\ce(p,q)=\se(p,q)\rtimes\bbR$ by $X_{ij}$, $X_i$ and~$X_0$.

\begin{rmk}\label{maxRemark}
{\rm 
It is worth noticing that the conformal Lie algebra $\so(p+1,q+1)$ is maximal in the
Lie algebra $\Vect_{\mathrm{Pol}}(\bbR^n)$ of polynomial vector fields in the
following sense: any bigger subalgebra of $\Vect_{\mathrm{Pol}}(\bbR^n)$
necessarily coincides with $\Vect_{\mathrm{Pol}}(\bbR^n)$. See \cite{BL} for a
simple proof. The uniqueness and the canonical character of our quantization procedure
definitely originates from this maximality property of $\so(p+1,q+1)$. 
See also \cite{Pos} for a classification of a class of maximal Lie subalgebras of
$\Vect_{\mathrm{Pol}}(\bbR^n)$.
}
\end{rmk}

\begin{rmk}\label{localRemark}
{\rm {}From now on, we will use local coordinate systems adapted to the flat
conformal structure on $M$ in which the generators of $\so(p+1,q+1)$ retain the form
(\ref{confGenerators}). 
This flat conformal structure precisely corresponds to a $\SO(p+1,q+1)$-structure on
$M$ (cf. Introduction) defined by the atlas of these adapted coordinate systems.
Clearly, our formul{\ae} will prove to be independent of the particular
choice of an adapted coordinate system and to be globally defined.  
}
\end{rmk}

\goodbreak

\subsection{Explicit formul{\ae} for the $\so(p+1,q+1)$-actions}\label{Comparison}

As a first application of the preceding results, let us compute the
action of the conformal algebra on $\cS_\d$ and on $\cD_{\l,\m}$ which is given by the
following two propositions.

\begin{pro}\label{confActionSymbPro}
The action of $\so(p+1,q+1)$ on $\cS_\d$ reads
\begin{equation}
\begin{array}{rcl}
\displaystyle L^\d_{X_i}  
&=&
\partial_i 
\\[10pt]
\displaystyle L^\d_{X_{ij}}
&=&
x_i\partial_j - x_j\partial_i
+
\xi_i\partial_{\xi^j} - \xi_j\partial_{\xi^i}
\\[10pt]
\displaystyle L^\d_{X_0}
&=&
x_i\partial_i
-
\xi_i\partial_{\xi_i}
+n\d
\\[10pt]
\displaystyle L^\d_{\bar{X}_i}
&=&
x_jx^j\partial_i - 2x_ix^j\partial_j
-2(\xi_i x_j - \xi_j x_i)\partial_{\xi_j}
+2\xi_j x^j\partial_{\xi^i}
-2n\d{}x_i.
\end{array}
\label{confActionSymbols}
\end{equation}
\end{pro}
\begin{proof}
These expressions follow from the explicit
form~(\ref{confGenerators}) of the $\so(p+1,q+1)$-generators, and from~(\ref{LieDer1}).
\end{proof}

\begin{pro}
The action of $\so(p+1,q+1)$ on $\cD_{\l,\m}$ reads
\begin{equation}
\cL_X^{\l,\m}
=
L_X^\d
\label{first}
\end{equation}
for all $X\in\ce(p,q)$, where $\d=\m-\l$; one furthermore has
\begin{equation}
\cL_{\bar X_i}^{\l,\m}=
L_{\bar X_i}^\d-
\xi_iT + 2(\cE+n\l)\,\partial_{\xi^i}
\label{InversionAction}
\end{equation}
for all infinitesimal inversions $\bar X_i$ (with $i=1,\ldots,n$) where
$T=\partial_{\xi^j}\partial_{\xi_j}$ is the trace and 
$\cE=\xi_j\partial_{\xi_j}$ the Euler operator. (See Section
\ref{Euclinv}.)
\end{pro}
\begin{proof}
This is a direct consequence of Proposition \ref{PolPro}, and of the
formul{\ae}~(\ref{confGenerators}).
\end{proof}

\begin{rmk}
{\rm 
The formula (\ref{InversionAction}) captures the difference between the 
$\so(p+1,q+1)$-modules $\cD_{\l,\m}$ and $\cS_\d$. Note that the operator
$\cL_X^{\l,\m}-L_X^\d$ is nilpotent since it maps $\cS^k_\d$ to $\cS^{k-1}_\d$.
}
\end{rmk}


\section{Main results}\label{MainResults}

In this section we formulate the main results of this article. All the proofs will be
given in Section \ref{ProofMain}.

\subsection{Quantization in the case $\delta=0$}

Let us consider first the special case $\d=\m-\l=0$ and use the shorthand notation
$\cS\equiv\cS_0$ and $\cD_\l\equiv\cD_{\l,\l}$. Although this is definitely not the
most general case to start with, this zero value of the shift is of central
importance to relate our conformally equivariant quantization to the more traditional
procedures such as geometric or deformation quantization.

\begin{thm}\label{ThmZeroShift}
(i) There exists an isomorphism of $\so(p+1,q+1)$-modules
\begin{equation}
\widetilde{\cQ}_\l:\cS\to\cD_\l.
\label{isomZeroShift}
\end{equation}
(ii) This isomorphism is unique provided the principal symbol be preserved at each
order, i.e., provided it reads $\widetilde{\cQ}_\l=\Id+\cN_\l$ with
nilpotent part $\cN_\l:\cS^k\to\cS^{k-1}$.
\end{thm}

Let us introduce a new operator on symbols that will eventually insure the
symmetry of the corresponding differential operators. 
Define
$\cI_\hbar:\cS_\d\to\cS_\d[i\hbar]$ by
\begin{equation}
\cI_\hbar(P)(\xi)=P(i\hbar\,\xi).
\label{Ihbar}
\end{equation}
Note that we will understand $\hbar$ either as a formal parameter or as a fixed real
number, depending upon the context. 

\begin{rmk}\label{quantumRemark}
{\rm
It is evident that $\cI_\hbar$ is an invariant operator, i.e., $[L^\d_X,\cI_\hbar]=0$
for all $X\in\Vect(\bbR^n)$. 
}
\end{rmk}

We then propose the
\begin{defi}
We will call conformally equivariant quantization the $\so(p+1,q+1)$-equivariant map
$\cQ_{\l;\hbar}:\cS\to\cD_\l[i\hbar]$ defined by
\begin{equation}
\cQ_{\l;\hbar}=\widetilde{\cQ}_\l\circ\cI_\hbar
\label{quantZeroShift}
\end{equation}
where $\hbar$ is a formal parameter and $\cI_\hbar$ is given by (\ref{Ihbar}).
\end{defi}

Theorem \ref{ThmZeroShift} and the preceding definition enable us to look for a
conformally invariant star-product on the space of symbols $\cS$ over $T^*M$. In
fact, as soon as one gets an isomorphism such as (\ref{quantZeroShift}), one can
readily define an associative bilinear operation (depending on $\l$)
\begin{equation}
*_{\l;\hbar}:\cS\otimes\cS\to\cS[[i\hbar]]
\label{starOperation}
\end{equation}
such that
\begin{equation}
\cQ_{\l;\hbar}(P*_{\l;\hbar}Q)
=
\cQ_{\l;\hbar}(P)\circ\cQ_{\l;\hbar}(Q).
\label{star}
\end{equation}

Recall that an associative operation $*_\hbar:\cS\otimes\cS\to\cS[[i\hbar]]$ is called
a star-product~\cite{BFFLS,dWL,Fed} if it is of the form
\begin{equation}
P*_\hbar{}Q=PQ+\frac{i\hbar}{2}\{P,Q\}+O(\hbar^2)
\label{starDef}
\end{equation}
where $\{\cdot,\cdot\}$ stands for the Poisson bracket on $T^*M$, and is given by
bi-differential operators at each order in~$\hbar$.

\begin{thm}\label{starProduct}
The associative, conformally invariant, operation $*_{\l;\hbar}$ defined by (\ref{star})
is a star-product if and only if $\l=\half$.
\end{thm}

Let us emphasize that this theorem provides us precisely with the value of
$\l$ used in geometric quantization and, in some sense, links the latter to
deformation quantization.

\subsection{General formulation. Resonant values of $\delta$}

In this section we formulate our result about the isomorphism of the
$\so(p+1,q+1)$-modules $\cS_\d$ and $\cD_{\l,\m}$ in the general situation.

The discussion below mainly relies on the structure of the spectrum of the Casimir
operator $\cC_{\l,\m}$ of (the $\so(p+1,q+1)$-module) $\cD_{\l,\m}$. Indeed, the
Casimir operator $\rC_\d$ of $\cS_\d$ turns out to be diagonalizable. Therefore, a
necessary condition for the $\so(p+1,q+1)$-modules $\cS_\d$ and $\cD_{\l,\m}$ to be
isomorphic is that $\cC_{\l,\m}$ be diagonalizable. This is of course the case if its
eigenvalues are ``simple'', while some problems could arise otherwise.
The latter case occurs only if the shift $\d=\m-\l$ belongs to the set
\begin{equation}
\Sigma=\left\{
\d_{k,\ell;s,t}\,\vert\,
k,\ell,s,t\in\bbN;\,
k>\ell;\,
2s\leq{}k;\,
2t\leq{}\ell
\right\}
\label{SigmaSet}
\end{equation}
where
\begin{equation}
\begin{array}{rcc}
\d_{k,\ell;s,t}
&=&
\displaystyle\frac{1}{n(k-\ell)}
\Big(
(k-\ell+t-s)(k+\ell-2(s+t)+n-1)\\[8pt]
&&+
(s-t)(k+\ell+1)
+
2(kt-\ell{}s)
\Big).
\label{deltaSpecial}
\end{array}
\end{equation}
The elements of $\Sigma$ will be called \textit{resonances}.

\begin{thm}(Generic case.)\label{IsomGen1}
If $n=p+q\geq2$ and $\d\not\in\Sigma$, then there exists an isomorphism of
$\so(p+1,q+1)$-modules
\begin{equation}
\widetilde{\cQ}_{\l,\m}:\cS_\d\to\cD_{\l,\m}
\label{prequant}
\end{equation}
which is unique provided the principal symbol be preserved at each order.
\end{thm}

If $\d\in\Sigma$, the Casimir operator $\cC_{\l,\m}$ has ``multiple'' eigenvalues
and, in some cases, is even not diagonalizable.
The corresponding \textit{critical} values of $\d$ are difficult to determine;
however, they belong to 
\begin{equation}
\Sigma_0=\left\{
\d_{k,\ell;s,t}\in\Sigma\,\vert\,
0\leq{}s-t\leq{}k-\ell
\right\}.
\label{SigmaSet0}
\end{equation}

\begin{thm}(Resonant case.)\label{IsomGen2}
If $n=p+q\geq2$ and $\d\in\Sigma\!\setminus\!\Sigma_0$,
then there exists an isomorphism (\ref{prequant}) of
$\so(p+1,q+1)$-modules
which is unique provided the principal symbol be preserved at each order.
\end{thm}

The proofs of Theorems \ref{IsomGen1} and \ref{IsomGen2} both consist mainly in
showing that $\cC_{\l,\m}$ is diagonalizable. In the case of Theorem \ref{IsomGen1},
this is quite immediate to prove, whereas the resonant case is much more involved.

\begin{rmk}
{\rm
If $\d\in\Sigma_0$, then there are values of the weights $\l$ and $\m$ for which the
sought isomorphism does exist (being, however, not necessarily unique). We have no
precise statement for this degenerate case, but in the example of second order
symbols, the table~(\ref{TheArray}) provides special values of 
$\l$ and $\m$ leading to an isomorphism~(\ref{prequant}).  
}
\end{rmk}

\begin{rmk}
{\rm
One easily finds values of $n$ for which $0\in\Sigma$ (for instance, $n=2$, for which
$\d_{4,3;2,0}=0$).
However, we will show that if $\d=0$ is resonant, it is not critical
(Lemma~\ref{deltaPositive}). }
\end{rmk}

\begin{rmk}
{\rm
In the one-dimensional case, $n=1$, the above theorems still hold
true but the resonances are simply $\d=1,\frac{3}{2},2,\frac{5}{2}$ and appear
in \cite{CMZ,Gar}.
}
\end{rmk}

Again, we will introduce the quantization map as the $\so(p+1,q+1)$-equivariant map
$\cQ_{\l,\m;\hbar}:\cS_\d\to\cD_\l[i\hbar]$ defined by
\begin{equation}
\cQ_{\l,\m;\hbar}=\widetilde{\cQ}_{\l,\m}\circ\cI_\hbar
\label{quantGen}
\end{equation}
as a natural generalization of (\ref{quantZeroShift}).

\goodbreak

Let us recall that if $\l+\m=1$, there exists, for compactly-supported
densities, a $\Vect(M)$-invariant pairing $\cF_\l\otimes\cF_\m\to\bbC$ defined by
\begin{equation}
\varphi\otimes\psi\mapsto\int_M{\!\overline{\varphi}\,\psi}.
\label{pairing}
\end{equation}

We can then formulate the important
\begin{cor}
Assume $\d\not\in\Sigma_0$ and $\l+\m=1$. The quantization
\begin{equation}
\check{P}
=
\cQ_{\frac{1-\d}{2},\frac{1+\d}{2};\hbar}(P)
\label{quantGenSymm}
\end{equation}
of any symbol $P\in\cS_\d$ is a symmetric (formally self-adjoint)
operator.
\end{cor}
\begin{proof}
Let us denote by $A^*$ the adjoint of $A\in\cD_{\l,1-\l}$ with respect to the
pairing~(\ref{pairing}). Consider the symmetric operator
$$
\Symm(\cQ_{\l,\m;\hbar}(P))
=
\half\Big(
\cQ_{\l,\m;\hbar}(P)+(\cQ_{\l,\m;\hbar}(P))^*
\Big)
$$
which exists whenever $\l+\m=1$. Notice that it has the same principal symbol as
$\cQ_{\l,\m;\hbar}(P)$. Now, the map $\Symm(\cQ_{\l,\m;\hbar})$ is obviously
$\so(p+1,q+1)$-equivariant. Theorems \ref{IsomGen1} and \ref{IsomGen2} just apply
and yield $\Symm(\cQ_{\l,\m;\hbar}(P))=\cQ_{\l,\m;\hbar}(P)$.
\end{proof}

The particular case $\d=0$ is of special importance and related to
Theorem~\ref{starProduct} since $\l=\m=\half$.

\subsection{Quantum Hamiltonians}\label{quantumHamiltonian}

To recover the traditional Schr\"odinger picture of quantum mechanics,
one needs to turn the operator $\check{P}$ resulting from our quantization
map (\ref{quantGenSymm}) into an operator 
\begin{equation}
\hat{P}:\cF_0\to\cF_0
\label{quantumOp}
\end{equation}
on the space of complex-valued functions on a conformally flat
manifold $(M,\rg)$.

\goodbreak

Using the natural identification $\cF_0\to\cF_\l$ between tensor densities and smooth
functions given (see Section \ref{VectModules}) by 
\begin{equation}
f\mapsto{}f\,\mod{\Vol_\rg}^\l,
\label{identification}
\end{equation}
one can introduce the differential operator, $\hat{P}$, defined by the commutative
diagram
\begin{equation}
\begin{CD}
\cF_0 @> \textstyle{\hat{P}} >> \cF_0 \strut\\
@V{\mod{\Vol_\rg}^\l}VV @VV{\mod{\Vol_\rg}^\m}V \strut\\
\cF_\l @> \textstyle{\check{P}} >> \cF_\m \strut
\end{CD}
\label{TheNewDiagram}
\end{equation}
where $\check{P}$ is given by (\ref{quantGenSymm}) in the case $\l+\m=1$. 

\begin{rmk}
{\rm 
So far, we only needed a conformal class of metrics to define a
conformally equivariant quantization map. But, in the current construction, we
definitely make a particular choice of metric, $\rg$, in the latter class to express the
operator $\hat{P}$.
}
\end{rmk}

In the case $\d=0$, which is relevant for quantum mechanics, the
operator~$\check{P}$ admits a prolongation as a (formally) self-adjoint operator on the
Hilbert space~$\overline{\cF_\half}$ (the completion of the space of compactly
supported half-densities with Hermitian inner product~(\ref{pairing})). It will be
therefore legitimate to call $\hat{P}$ the quantum Hamiltonian associated with the
Hamiltonian $P\in\cS\cong\Pol(T^*\bbR^n)$. This quantum Hamiltonian is then a (formally)
self-adjoint operator on the space
$L^2(M,|\Vol_\rg|)$.

\section{Conformally invariant operators}\label{ConfInvOp}

The space $\bbC[x^1,\ldots,x^n,\xi_1,\ldots,\xi_n]$ of polynomials on $T^*\bbR^n$ is
naturally a module over the Lie algebra, $\Vect_{\mathrm{Pol}}(\bbR^n)$, of polynomial
vector fields on $\bbR^n$. This module structure is induced by the $\Vect(\bbR^n)$
action on $T^*\bbR^n$. But, we will rather consider, as in Section~\ref{TheModules}, the
deformed action (\ref{LieDer1}) depending on a parameter $\d$; we will henceforth
denote this module by $\bbC[x^1,\ldots,x^n,\xi_1,\ldots,\xi_n]_\d$.

\begin{defi}
We denote by
$\End_\diff(\bbC[x^1,\ldots,x^n,\xi_1,\ldots,\xi_n])$ the subspace
$$
\bbC[
x^1,\ldots,x^n,
\xi_1,\ldots,\xi_n,
\frac{\partial}{\partial{}x^1},\ldots,\frac{\partial}{\partial{}x^n},
\frac{\partial}{\partial{}\xi_1},\ldots,\frac{\partial}{\partial{}\xi_n}
]
$$
of polynomial differential operators on 
$\bbC[x^1,\ldots,x^n,\xi_1,\ldots,\xi_n]$.
\end{defi}

\begin{defi}
To any Lie algebra $\fg\subset\Vect_{\mathrm{Pol}}(\bbR^n)$ we associate its
com\-mutant,~$\fg^!$, as the Lie subalgebra of
$\End_\diff(\bbC[x^1,\ldots,x^n,\xi_1,\ldots,\xi_n])$ of those operators that commute with
$\fg$. 
\end{defi}
This classical notion of commutant has first been considered in the context of
differential operators by Kirillov \cite{Kir2}.

\subsection{Algebra of Euclidean invariants}\label{Euclinv}

To work out a conformally equivariant quantization map, we need to study first
equivariance with respect to the Euclidean subalgebra $\se(p,q)$. To this end, we will
introduce the commutant $\se(p,q)^!$.

Let us recall the structure of $\se(p,q)^!$ which has been shown
\cite{DO2} to be the associative algebra generated by the operators
\begin{equation}
\rR=\xi^i\xi_i,\qquad
\rE=\xi_i\frac{\partial}{\partial \xi_i}+\frac{n}{2},\qquad
\rT=\frac{\partial}{\partial \xi^i}\frac{\partial}{\partial \xi_i}
\label{sl2}
\end{equation}
whose commutation relations are those of $\Sl(2,\bbR)$ together with
\begin{equation}
\rG=\xi^i\frac{\partial}{\partial x^i},\qquad
\rD=\frac{\partial}{\partial \xi_i}\frac{\partial}{\partial x^i},\qquad
\Delta=\frac{\partial}{\partial x^i}\frac{\partial}{\partial x_i}
\label{rh1}
\end{equation}
which generate the Heisenberg Lie algebra $\rh_1$.

\goodbreak

We will find it useful to deal with the Euler operator
\begin{equation}
\Euler=\rE-\frac{n}{2}.
\label{Euler}
\end{equation}

An example of Howe dual pairs of (non semi-simple) Lie algebras is given by
\begin{thm}\cite{DO2}
\label{Usp3}
The commutant $\se(p,q)^!$ in $\End_\diff(\bbC[x^1,\ldots,x^n,\xi_1,\ldots,\xi_n])$ is
isomorphic to 
$U({\Sl(2,\bbR)\ltimes\rh_1})/\cZ$ where the ideal, $\cZ$, is as follows
\hfill\break
(i) if $n=2$, the ideal $\cZ$ is generated by
\begin{equation}
Z =
\Big(\rC+{
\frac{3}{2}
}\Big)\,\Delta
+
{
\frac{1}{4}
}\Big(\left[\rD,[\rG,\rC]\right]_+
-
\left[\rG,[\rD,\rC]\right]_+
\Big),
\label{Rel2}
\end{equation}
where $[\,\cdot\,,\cdot\,]_+$ stands for the anticommutator, and
\begin{equation}
\rC=\rE^2-\frac{1}{2}[\rR,\rT]_+
\label{CasimirSl2}
\end{equation}
is the Casimir element of $\Sl(2,\bbR)$,
\hfill\break
(ii) if $n\geq3$, one has
\begin{equation}
\cZ = \{0\}.
\label{Rel3}
\end{equation}
\end{thm}

This theorem is a generalization of the celebrated Brauer-Weyl Theorem \cite{Wey} (see
also \cite{Kir2,DO2}).

Let us mention that we will, actually, need considering invariant operators with
respect to homotheties generated by $X_0$ (see~(\ref{confGenerators}) and
(\ref{subalgebras})) inside $\se(p,q)^!$. We readily have the 
\begin{cor}\label{corComm}
The commutant $\ce(p,q)^!$ in $\End_\diff(\bbC[x^1,\ldots,x^n,\xi_1,\ldots,\xi_n]_\d)$ is the
associative algebra generated (see (\ref{sl2}) and (\ref{rh1})) by
\begin{equation}
\rE,\qquad 
\rR_0=\rR\,\rT,\qquad 
\rD,\qquad
\rG_0=\rG\,\rT,\qquad
\Delta_0=\Delta\,\rT.
\label{homInv}
\end{equation}
\end{cor}

\subsection{Algebra of conformal invariants}\label{confinv}

The commutant $\so(p+1,q+1)^!$, is given by the following corollary of
Theorem \ref{Usp3}.
\begin{cor}\label{noinvariants}
The commutant $\so(p+1,q+1)^!$ in $\End_\diff(\bbC[x^1,\ldots,x^n,\xi_1,\ldots,\xi_n]_\d)$ is,
for $n\geq3$, the commutative associative algebra generated by $\rE$ and $\rR_0$.
\end{cor}
\begin{proof}
In view of the preceding corollary, we need only the commutation
relations of the operators (\ref{homInv}) with the generators $L^\d_{\bar{X}_i}$ of
inversions given in (\ref{confActionSymbols}) in order to determine
$\so(p+1,q+1)^!$.
\goodbreak
Straightforward calculation leads to
\begin{eqnarray}
[\rE,L^\d_{\bar{X}_i}] &=& 0\nonumber
\\[6pt]
[\rR_0,L^\d_{\bar{X}_i}]  &=& 0\nonumber
\\[6pt]
[\rG_0,L^\d_{\bar{X}_i}]  
&=& 
2\Big(\rR_0\partial_{\xi^i} + (2-n\d)\xi_i\rT\Big)\label{G0}
\\[6pt]
[\rD,L^\d_{\bar{X}_i}] 
&=& 
2\Big(-\xi_i\rT +2\cE\partial_{\xi^i} + n(1-\d)\partial_{\xi^i}\Big)\nonumber
\\[6pt]
[\Delta_0,L^\d_{\bar{X}_i}]  
&=& 
4\Big(
\cE\partial_i\rT + \rG_0\partial_{\xi^i} - \xi_i\rD\rT\Big) 
+ 2\Big(2+n(1-2\d)\Big)\partial_i\rT
\nonumber
\label{CommutationRelations}
\end{eqnarray}
for $i=1,\ldots,n$.

Now, Theorem \ref{Usp3} guarantees that the monomials
$\rW_{e,r,d,g,\ell}=\rE^e\,\rR_0^r\,\rD^d\,\rG_0^g\,\Delta_0^\ell$ are independent;
they are of degree 
\begin{equation}
\deg(\rW_{e,r,d,g,\ell})=e+2r+d+2g+2\ell
\label{degree}
\end{equation}
as differential operators in $\xi_i$. We thus have to look for the commutant of
inversions in $\ce(p,q)^!$, i.e., to determine the linear combinations of the
previous monomials that commute with $L^\d_{\bar{X}_i}$. The commutator
$[\rW_{e,r,d,g,\ell},L^\d_{\bar{X}_i}]$ is a differential operator (in $\xi_i$) of
degree $\deg(\rW_{e,r,d,g,\ell})+1$. We are therefore led to study the principal symbol
of this operator, which can be easily computed with the help of (\ref{G0}).
In order to make our calculations more tractable, let us rather deal with the principal
symbol of the operator $\sum_{i=1}^n{\xi^i[\rW_{e,r,d,g,\ell},L^\d_{\bar{X}_i}]}$; it is
of the form
\begin{equation}
\begin{array}{l}
2g\,\rE^{e+1}\,\rR_0^{r+1}\,\rD^d\,\rG_0^{g-1}\,\Delta_0^\ell
\\[8pt]
-2d\Big(\rE^e\,\rR_0^{r+1}\,\rD^{d-1}\,\rG_0^g\,\Delta_0^\ell
-2\,\rE^{e+2}\,\rR_0^r\,\rD^{d-1}\,\rG_0^g\,\Delta_0^\ell
\Big)
\\[8pt]
+4\ell\Big(
2\,\rE^{e+1}\,\rR_0^r\,\rD^d\,\rG_0^{g+1}\,\Delta_0^{\ell-1}
-\rE^e\,\rR_0^{r+1}\,\rD^{d+1}\,\rG_0^g\,\Delta_0^{\ell-1}
\Big).
\end{array}
\label{PrincipalSymbol}
\end{equation}
We then seek the linear combinations of the monomials~$\rW_{e,r,d,g,\ell}$, of fixed
degree (\ref{degree}), for which the previous expression is identically zero.
Resorting to Theorem~\ref{Usp3}, we immediately get $g=0$ since the first term in
(\ref{PrincipalSymbol}) is clearly independent of the others. The same is true for the
next two terms, yielding $d=0$ and $\ell=0$.
\end{proof}

\goodbreak

\subsection{Casimir operator $\rC_\d$}

We have computed in Section \ref{confinv} the commutant $\so(p+1,q+1)^!$ of the
conformal Lie algebra. Now, representation theory tells us that there exists a
distinguished invariant within this commutant, namely the Casimir operator.
 
Recall that the Casimir operator of a given representation
$\rho:\fg\to\End(V)$ of a semi-simple Lie algebra $\fg$ is
\begin{equation}
\rC=B^{\a\b}\rho({X_\a})\rho({X_\b})
\end{equation}
where $B$ is the Killing metric and $(X_\a)$
any basis of $\fg$ with $B_{\a\b}=B(X_\a,X_\b)$, the components of the associated
Gram matrix.
It is well known that the Casimir operator is invariant, i.e.
$[\rC,\rho(X)]=0$ for all $X\in\fg$.

In this section, we will provide the explicit calculation of the Casimir operator
of the $\so(p+1,q+1)$-action, $L^\d$, on $\cS_\d$ given by Proposition
\ref{confActionSymbPro}. We choose the Killing form as 
\begin{equation}
B(X,Y)=-\frac{1}{2}\Tr(XY)
\label{Killing}
\end{equation}
where $X,Y\in\so(p+1,q+1)$ in their $(n+2)\times(n+2)$ matrix realization. 
We can then give the explicit formula for this Casimir operator, denoted by
$\rC_\d$, in terms of the invariant operators (\ref{homInv}). 
\begin{pro}\label{casimirPro}
One has
\begin{equation}
\rC_\d=\rR_0+2(1+n(\d-1)-\cE)\cE-n^2\d(\d-1)
\label{casimir}
\end{equation}
where $n=p+q$.
\end{pro}
\goodbreak
\begin{proof}
The matrix realization of the $\so(p+1,q+1)$-generators
(\ref{confGenerators}) is given by
\begin{eqnarray}
X_i
=
\pmatrix{
0&-\sqrt{2}e_i&0\cr
0&0&0\cr
\sqrt{2}e^\flat_i&0&0\cr
},
&&\\[10pt]
X_{ij}
=
\pmatrix{
e_je^\flat_i-e_ie^\flat_j&0&0\cr
0&0&0\cr
0&0&0\cr
},
&&
X_0
=
\pmatrix{
0&0&0\cr
0&-1&0\cr
0&0&1\cr
},\\[10pt]
\bar{X}_i
=
\pmatrix{
0&0&\sqrt{2}e_i\cr
-\sqrt{2}e^\flat_i&0&0\cr
0&0&0\cr
},
&&
\end{eqnarray}
where $(e_i)$ is the canonical basis of $\bbR^n$ and $(e^\flat_i=g(e_i))$ its dual
basis associated with the metric $g$.

A simple calculation yields the basis $(X^\a=B^{\a\b}X_\b)$ of $\so(p+1,q+1)$ dual
to~$(X_\a)$ with respect to the Killing metric (\ref{Killing}). One gets
\begin{equation}
\begin{array}{rcl}
X^i   
&=&
-\frac{1}{2}\,g^{ij}\bar{X}_j \\[10pt]
X^{ij}
&=&
g^{ik}g^{j\ell}X_{k\ell} \\[10pt]
X^0
&=&
-X_0 \\[10pt]
\bar{X}^i
&=&
-\frac{1}{2}\,g^{ij}X_j.
\end{array}
\label{dualconfGenerators}
\end{equation}
Using the $\so(p+1,q+1)$-action, $L^\d$, on $\cS_\d$ given in
(\ref{confActionSymbols}), one shows immediately that the Casimir operator
\begin{equation}
\rC_\d
=
\half{}g^{ik}g^{j\ell}L^\d_{X_{ij}}L^\d_{X_{k\ell}}
-\left(L^\d_{X_0}\right)^2
-\half{}g^{ij}L^\d_{X_i}L^\d_{\bar{X}_j}
-\half{}g^{ij}L^\d_{\bar{X}_i}L^\d_{X_j}
\label{explicitcasimir}
\end{equation}
actually retains the form (\ref{casimir}).
\end{proof}

\goodbreak

\begin{rmk}\label{CasimirRemark}
{\rm 
It is worth noticing that the Casimir operator (\ref{casimir}) can be alternatively
expressed in terms of the Casimir operator, $\rC$, (see (\ref{CasimirSl2})) and the
Cartan generator, $\rE$, of $\Sl(2,\bbR)$. One finds
\begin{equation}
\rC_\d=-\rC-(\rE-n\d)^2-n\left(1-\frac{n}{2}\right).
\end{equation}
}
\end{rmk}

\subsection{Casimir operator $\cC_{\l,\m}$}\label{CasimirSection}

The Casimir operator of the $\so(p+1,q+1)$-action on $\cD_{\l,\m}$ is defined,
accordingly, by
\begin{equation}
\cC_{\l,\m}=B^{\a\b}\cL^{\l,\m}_{X_\a}\cL^{\l,\m}_{X_\b}.
\label{defCasimir}
\end{equation}

\begin{pro}\label{CasimirPro}
The Casimir operator of the $\so(p+1,q+1)$-action on $\cD_{\l,\m}$ is of the form
\begin{equation}
\cC_{\l,\m}=\rC_\d + \rG_0 - 2(n\l+\cE)\rD.
\label{Casimir}
\end{equation}
\end{pro}
\begin{proof}
The explicit formula for $\cC_{\l,\m}$ is obviously obtained by replacing $L^\d$ by
$\cL^{\l,\m}$ in (\ref{explicitcasimir}). Applying then (\ref{first}) and
(\ref{InversionAction}) to that expression immediately leads to the sought result
(\ref{Casimir}).
\end{proof}

\section{Proofs of the Main Theorems}\label{ProofMain}

Throughout this section we freely use, for convenience, the local
identification~(\ref{sigma}) of $\cD_{\l,\m}$ and $\cS_\d$.

\subsection{Diagonalization of the Casimir operator $\rC_\d$}

We have already mentioned that we will study the diagonalization of the Casimir
operators $\rC_\d$ and $\cC_{\l,\m}$. Here, we understand that an
endomorphism of an infinite-dimensional space is diagonalizable if any element of
the latter is a (finite) sum of eigenvectors of the former.

Let us recall that $[\cE,\rR_0]=0$, so that $[\rC_\d,\cE]=0$ and $[\rC_\d,\rR_0]=0$.
We can thus simultaneously diagonalize the three operators $\cE,\rR_0$ and
$\rC_\d$.

\begin{lem}\label{eigenR0}
The eigenvectors of the operator $\rR_0$ restricted to the space $\cS_{k,\d}$ of
homogeneous polynomials (see (\ref{grad})) are of the form
\begin{equation}
P_{k,s} = \rR^s Q,
\label{Pks}
\end{equation}
where $\rR$ is given by (\ref{sl2}) and $Q\in\cS_{k-2s,\d}$
is tracefree (harmonic), viz $\rT{}Q=0$, and 
\begin{equation}
s\in\left\{0,1,2,\ldots,\left[k/2\right]\right\}.
\label{sRange}
\end{equation}
The associated eigenvalues are
\begin{equation}
\varrho_{k,s}=2s(n+2(k-s-1)).
\label{rhoks}
\end{equation}
\end{lem}
\noindent
For a proof, see, e.g., \cite{Wey}. 

\medskip
\goodbreak

We readily have the following
\begin{cor}
The spectrum of the Casimir operator $\rC_\d$ is given by
\begin{equation}
\gamma_{k,s}=2s(n+2(k-s-1)) + 2k(1+n(\d-1)-k) -n^2\d(\d-1).
\label{gammaks}
\end{equation}
\end{cor}
\begin{proof}
This follows immediately from Lemma \ref{eigenR0} and the expression
(\ref{casimir}) of the Casimir operator~$\rC_\d$.
\end{proof}

\goodbreak

We have thus the following useful decomposition: every $P\in\cS_\d$ can be
decomposed as a locally finite sum
\begin{equation}
P=\sum_{k=0}^\infty{\sum_{s=0}^{[k/2]}{P_{k,s}}}
\label{decomp}
\end{equation}
where $[k/2]$ is the integer part of $k/2$.
In other words, we have a direct sum decomposition
\begin{equation}
\cS_\d=
\bigoplus_{\stackrel{\scriptstyle{}k=0}{s\leq[k/2]}}^\infty{\cS_{(k,s),\d}}
\label{DirSum}
\end{equation}
into eigenspaces of $\cE$ and $\rR_0$. (See Theorem (5.6.A) in \cite{Wey}.)

\medskip

We are now able to explain the origin of the resonant values (\ref{deltaSpecial})
of~$\d=\m-\l$.

\begin{lem}\label{Obtainingdelta}
One has $\d=\d_{k,\ell;s,t}$
if and only if $\gamma_{k,s}=\gamma_{\ell,t}$.
\end{lem} 
\begin{proof}
This is straightforward from (\ref{gammaks}).
\end{proof}

\subsection{Diagonalization of the Casimir operator $\cC_{\l,\m}$}

We establish, in this section, the main technical statement that helps us to prove
the existence of an equivariant quantization map for almost all
values of the shift~$\d$. 


The expression (\ref{Casimir}) of the Casimir operator is of the form 
$\cC_{\l,\m}=\rC_\d + \rN_\l$ with nilpotent part $\rN_\l:\cS^k_\d\to\cS^{k-1}_\d$
(for $k=0,1,2,\ldots$),  see Section \ref{CasimirSection}. 
This implies that any solution $P\in\cD^k_{\l,\m}$ of
the equation
\begin{equation}
\cC_{\l,\m}P=\gamma P
\label{theEigenEquation}
\end{equation}
is of the form  $P=P_{k,s}+P'$
where
$P_{k,s}$ is as in (\ref{Pks}) and
$P'\in\cD^{k-1}_{\l,\m}$; the eigenvalue~$\gamma$ clearly coincides
with $\gamma_{k,s}$ given by (\ref{gammaks}).

\begin{pro}(Generic case.) \label{DiagCasimir1}
If $\d\not\in\Sigma$, then the eigenvalue problem
(\ref{theEigenEquation}) for the Casimir operator $\cC_{\l,\m}$ has a
solution if and only if $\gamma=\gamma_{k,s}$ for some $k,s$ (as given by
(\ref{gammaks})). The corresponding eigenvectors are uniquely determined by their
principal symbols, arbitrarily taken in $\cS_{(k,s),\d}$.
\end{pro}
\begin{proof}
The highest degree component of the eigenvalue equation (\ref{theEigenEquation}) is
just the eigenvalue equation for $\rC_\d$.  Hence, a solution of
(\ref{theEigenEquation}) is necessarily of the form $\gamma=\gamma_{k,s}$ and 
\begin{equation}
P=P_{k,s}+
\sum_{\stackrel{\scriptstyle\ell<k}{t\leq[\ell/2]}}{P_{\ell,t}}
\label{theDecomposition}
\end{equation}
according to the decomposition (\ref{DirSum}).
The remainder of equation (\ref{theEigenEquation}) reads now
\begin{equation}
\sum_{\stackrel{\scriptstyle\ell<k}{t\leq[\ell/2]}}{
(\g_{k,s}-\g_{\ell,t})P_{\ell,t}}
=
\rN_\l P.
\label{theRecurrence}
\end{equation}
Since $\d\not\in\Sigma$, by Lemma \ref{Obtainingdelta}, the coefficients
$\g_{k,s}-\g_{\ell,t}$ do not vanish. In view of the nilpotency of~$\rN_\l$, the
result follows immediately.
\end{proof}

In order to handle the case of non-critical resonant values of $\d$, we need the
following
\begin{lem}\label{DecreasingDegree}
If $P\in\cS_{(k,s),\d}$, then 

\noindent
(i) the polynomials $\rD(P)$ and $\rG_0(P)$ belong to
$\cS_{(k-1,s-1),\d}\oplus\cS_{(k-1,s),\d}$,

\noindent
(ii) the polynomial $\Delta_0(P)$ belongs to
$\cS_{(k-2,s-2),\d}\oplus\cS_{(k-2,s-1),\d}\oplus\cS_{(k-2,s),\d}$.
\end{lem} 
\begin{proof}
Any polynomial $P\in\cS_{(k,s),\d}$ is of the form $P=\rR^sQ$ where
$Q\in\cS_{(k-2s,0),\d}$ is harmonic, see (\ref{Pks}). One has
$$
\begin{array}{rcl}
\rD(P)
&=&
\displaystyle
[\rD,\rR^s](Q)+\rR^s\rD(Q)\\[6pt]
&=&
\displaystyle\sum_{r=0}^{s-1}{\rR^r[\rD,\rR]\rR^{s-r-1}(Q)}
+\rR^s\rD(Q)\\[16pt]
&=&
\displaystyle{}2s\rR^{s-1}\rG(Q)+\rR^s\rD(Q)
\end{array}
$$
since $[\rD,\rR]=2\rG$ and $[\rR,\rG]=0$ in the Lie algebra
$\Sl(2,\bbR)\ltimes\rh_1$ (see Section \ref{Euclinv}).

At last, $\rD(Q)$ is harmonic because $[\rT,\rD]=0$ and one furthermore easily checks
that 
$\rG(Q)\in\cS_{(k-2s+1,0),\d}\oplus\cS_{(k-2s+1,1),\d}$. 
Hence,
$\rD(P)\in\cS_{(k-1,s-1),\d}\oplus\cS_{(k-1,s),\d}$.

The proof for $\rG_0(P)$ and $\Delta_0(P)$ is analogous and will be omitted.
\end{proof}

Let us then introduce the space
\begin{equation}
\widetilde{\cS}_{(k,s),\d}
=
\bigoplus_{0\leq{}s-t\leq{}k-\ell}{\cS_{(\ell,t),\d}}
\end{equation}
which is ``generated'' by the tree
\begin{equation}
\begin{array}{rcl}
& \cS_{(k,s),\d} & \\
& \swarrow\searrow &\\
\cS_{(k-1,s-1),\d} && \cS_{(k-1,s),\d}\\
\swarrow \searrow &&\swarrow \searrow\\
\ldots \hfill \quad &\cS_{(k-2,s-1),\d}& \hfill \ldots\\
&\ldots&
\end{array}
\label{tree}
\end{equation}
In view of the preceding lemma, $\widetilde{\cS}_{(k,s),\d}$ is stabilized by the
Casimir operator $\cC_{\l,\m}$. Moreover, if $\d\not\in\Sigma_0$ (see
(\ref{SigmaSet0})), $\g_{\ell,t}\neq\g_{k,s}$ whenever
$\cS_{(\ell,t),\d}\subset\widetilde{\cS}_{(k,s),\d}$.

\begin{pro}(Resonant case.) \label{DiagCasimir2}
If $\d\in\Sigma\!\setminus\!\Sigma_0$, 
then the eigenvalue problem
(\ref{theEigenEquation}) for the Casimir operator $\cC_{\l,\m}$ has a
solution if and only if $\gamma=\gamma_{k,s}$ for some $k,s$ (as given by
(\ref{gammaks})). 
Every $P_{k,s}\in\cS_{(k,s),\d}$ is the principal
symbol of an eigenvector of~$\cC_{\l,\m}$. This eigenvector is
uniquely determined provided it belongs to $\widetilde{\cS}_{(k,s),\d}$.
\end{pro}
\begin{proof}
One proceeds just as in the proof of Proposition \ref{DiagCasimir1}.
The fact that one restricts $\ell,t$ in the decomposition 
(\ref{theDecomposition}) to those values for which
$\cS_{(\ell,t),\d}\subset\widetilde{\cS}_{(k,s),\d}$ again insures that the
coefficients in equation (\ref{theRecurrence}) do not vanish.
\end{proof}

\begin{rmk}\label{ResonantRemark}
{\rm 
In Proposition \ref{DiagCasimir2}, if $\d\neq\d_{k,\ell;s,t}$ for some $\ell,t$, then
any eigenvector with principal symbol in $\cS_{(k,s),\d}$ necessarily belongs to
$\widetilde{\cS}_{(k,s),\d}$.  }
\end{rmk}

If $\d$ is not critical, it then follows from Propositions \ref{DiagCasimir1} and
\ref{DiagCasimir2} that  every $P_{k,s}\in\cS_{(k,s),\d}$ is the principal symbol of
an eigenvector $\widetilde{P}_{k,s}\in\widetilde{\cS}_{(k,s),\d}$ of the Casimir
operator
$\cC_{\l,\m}$. Hence the 
\begin{cor}
If $\d\not\in\Sigma_0$ the Casimir operator $\cC_{\l,\m}$ is diagonalizable.
\end{cor}

\goodbreak

\subsection{Proof of Theorems \ref{IsomGen1} and \ref{IsomGen2}}

Let us show that the diagonalization of the Casimir operator actually leads to the
determination of a unique isomorphism of
$\so(p+1,q+1)$-modules
$\widetilde{\cQ}_{\l,\m}:\cS_\d\to\cD_{\l,\m}$. Looking for a map
$\widetilde{\cQ}_{\l,\m}$ such that
\begin{equation}
\begin{CD}
\cD_{\l,\m} @> \textstyle{\cC_{\l,\m}} >> \cD_{\l,\m} \strut\\
@A{\textstyle{\widetilde{\cQ}_{\l,\m}}}AA 
@AA{\textstyle{\widetilde{\cQ}_{\l,\m}}}A \strut\\
\cS_\d @> \textstyle{\rC_\d} >> \cS_\d \strut
\end{CD}
\label{TheDiagram}
\end{equation}
be a commutative diagram, we are led to the
\begin{defi}
The linear map $\widetilde{\cQ}_{\l,\m}:\cS_\d\to\cD_{\l,\m}$ is defined
by 
\begin{equation}
\widetilde{\cQ}_{\l,\m}(P_{k,s})=\widetilde{P}_{k,s}
\label{Quant}
\end{equation}
using the decomposition (\ref{decomp}).
\end{defi}
This map has, obviously, the following properties
\begin{enumerate}
\item
$\cC_{\l,\m}\widetilde{\cQ}_{\l,\m}=\widetilde{\cQ}_{\l,\m}\rC_\d$,\label{a}
\item
$\widetilde{\cQ}_{\l,\m}=\Id+\cN_{\l,\m}$ with nilpotent part
$\cN_{\l,\m}:\cS^k_\d\to\cS^{k-1}_\d$\label{b}.
\end{enumerate}

\subsubsection{Proof of Theorem \ref{IsomGen1}}

Let us first prove that
\begin{equation}
\cL^{\l,\m}_X\widetilde{\cQ}_{\l,\m}=\widetilde{\cQ}_{\l,\m}L^\d_X
\label{EquivProp}
\end{equation}
for all $X\in\so(p+1,q+1)$.
{}From Property \ref{a}, we see that if $P_{k,s}\in\cS_{(k,s),\d}$ 
then $\cL^{\l,\m}_X\widetilde{\cQ}_{\l,\m}P_{k,s}$ and
$\widetilde{\cQ}_{\l,\m}L^\d_XP_{k,s}$ are both eigenvectors of the Casimir operator
$\cC_{\l,\m}$ associated with the same eigenvalue $\gamma_{k,s}$. Moreover,
Property
\ref{b} and (\ref{LieDer2}) entail that
$\cL^{\l,\m}_X\widetilde{\cQ}_{\l,\m}P_{k,s}$ and
$\widetilde{\cQ}_{\l,\m}L^\d_XP_{k,s}$ have the same principal symbol, namely
$L^\d_XP_{k,s}$. It follows from Proposition
\ref{DiagCasimir1} that these eigenvectors actually coincide. Hence, the existence
of the sought quantization map. 

Now, to prove the uniqueness, it suffices to note that an isomorphism of the
$\so(p+1,q+1)$-modules $\cS_\d$ and $\cD_{\l,\m}$ necessarily intertwines the
corresponding Casimir operators. If it moreover preserves the principal symbol, then
Proposition
\ref{DiagCasimir1} shows that it is, indeed,
$\widetilde{\cQ}_{\l,\m}$.

The proof of Theorem \ref{IsomGen1} is complete.

\subsubsection{Proof of Theorem \ref{IsomGen2}}

This proof is built on the same pattern as the previous one.
But, since $\d$ has resonant values, we must resort to Proposition
\ref{DiagCasimir2} instead of Proposition \ref{DiagCasimir1}. This is done at the
expense of some preparation due to the fact that the uniqueness of the eigenvector
$\widetilde{P}_{k,s}$ is guaranteed only within $\widetilde{\cS}_{(k,s),\d}$.

The above proof of the equivariance property (\ref{EquivProp}) should now be
completed with the help of
\begin{lem}
For every $X\in\so(p+1,q+1)$, one has $L^\d_X\cS_{(k,s),\d}\subset\cS_{(k,s),\d}$
and $\cL^{\l,\m}_X\widetilde{\cS}_{(k,s),\d}\subset\widetilde{\cS}_{(k,s),\d}$.
\end{lem}
\begin{proof}
The first inclusion easily follows from Corollary \ref{noinvariants}. As to the
second one, we then proceed as in the proof of Lemma
\ref{DecreasingDegree}, using (\ref{InversionAction}).
\end{proof}
The existence of the isomorphism $\widetilde{\cQ}_{\l,\m}$ is thus proven.

In the same way, the uniqueness of the sought isomorphism is established as in the
proof of Theorem \ref{IsomGen1}, provided we apply the following
\begin{lem}
Any linear map $\widetilde{\cQ}:\cS_\d\to\cS_\d$ that intertwines the
$\ce(p,q)$-action and do not increase the degree stabilizes each space
$\widetilde{\cS}_{(k,s),\d}$.
\end{lem}
\begin{proof}
It has been shown in Lemma 7.1 of \cite{Lec1} (see also Theorem 5.1 of \cite{LO})
that such a map $\widetilde{\cQ}$ is necessarily a differential operator with
constant coefficients.
We can thus apply Corollary \ref{corComm}: $\widetilde{\cQ}$ is a polynomial in the
operators (\ref{homInv}). We conclude by using Lemma \ref{DecreasingDegree}.
\end{proof}

This ends the proof of Theorem \ref{IsomGen2}.

\subsection{Proof of Theorem \ref{ThmZeroShift}}

In order to prove Theorem \ref{ThmZeroShift}, it is enough to show that $\d=0$ is
not a critical value, i.e., $0\not\in\Sigma_0$. This follows from the stronger
\begin{lem}\label{deltaPositive}
If the following inequalities hold
\begin{equation}
0\leq{}s-t\leq{}k-\ell,
\label{bad}
\end{equation}
one has $\d_{k,\ell;s,t}>0$.
\end{lem} 
\begin{proof}
In the expression (\ref{deltaSpecial}), both factors in the first
term, $k-\ell+t-s$ and $k+\ell-2(s+t)+n-1$, are non-negative in view of (\ref{bad}) and
(\ref{sRange}). Also (\ref{bad}) yields $kt-\ell{}s\geq-\ell(s-t)$ so that the second
term is bounded from below by
$(s-t)(k-\ell+1)$, which is non-negative.  We have just shown that $\d\geq0$. 

Now, if $s=t$, one has $n\d=k+\ell+n-2t-1\geq{}k$ because of (\ref{sRange}). The
result follows since $k>0$ in (\ref{deltaSpecial}).
\end{proof}

\subsection{Proof of Theorem \ref{starProduct}: \\Conformally invariant star-\-product}

Throughout this section we will only consider the case $\l=\m$. Let us give an explicit
expression for the quantization map (\ref{quantZeroShift}) up to the second order in
$\hbar$.

\goodbreak

\begin{pro}
\label{firstOrderIsom}
If $P_{k,s}\in\cS_{(k,s)}$ is a homogeneous polynomial (see (\ref{DirSum})) with $k>2$,
then

\noindent
(i) if $s>0$, the quantization map is of the form
\begin{equation}
\cQ_{\l;\hbar}(P_{k,s})=
P_{k,s}+\frac{i\hbar}{2}\left(
\rD(P_{k,s})+
\frac{(1-2\l)n}{s(2s-2k-n+2)}\rG_0(P_{k,s})
\right)+O(\hbar^2),
\label{firstOrderIsomForm}
\end{equation}

\noindent
(ii) if $s=0$ (i.e., the harmonic case), one has
\begin{equation}
\cQ_{\l;\hbar}(P_{k,0})=
P_{k,0}+i\hbar\left(\frac{n\l+k-1}{n+2(k-1)}\rD(P_{k,0})\right)+
O(\hbar^2).
\label{firstOrderIsomFormHarm}
\end{equation}
\end{pro}
\begin{proof}
(i) An eigenvector of the Casimir operator $\cC_{\l,\l}$, with principal symbol
$P_{k,s}$, is of the form
$$
P=P_{k,s}+P_{k-1,s}+P_{k-1,s-1}+\hbox{terms of degree\ }\leq{k-2}.
$$
(See the formula (\ref{Casimir}) and Lemma \ref{DecreasingDegree}.)

The eigenvalue
problem (\ref{theEigenEquation}) therefore leads to
$\g_{k,s}P_{k,s}=\g{}P_{k,s}$, where $\g_{k,s}$ is the eigenvalue of the Casimir
operator $\rC_0$ given by (\ref{gammaks}), and to
\begin{equation}
\begin{array}{rcl}
\g_{k-1,s}P_{k-1,s}+\g_{k-1,s-1}P_{k-1,s-1}
&=&
\g(P_{k-1,s}+P_{k-1,s-1})\\[10pt]
&&+
2(n\l+k-1)\rD(P_{k,s})-\rG_0(P_{k,s}).
\end{array}
\label{NextLevel}
\end{equation}
In order to solve this equation for $P_{k-1,s}$ and $P_{k-1,s-1}$, one needs
to introduce the projectors
$$
\Pi^{k-1,s}_1=
\frac{\rR_0-\vr_{k-1,s-1}}{\vr_{k-1,s}-\vr_{k-1,s-1}}
\qquad\hbox{and}\qquad
\Pi^{k-1,s}_2=
\frac{\rR_0-\vr_{k-1,s}}{\vr_{k-1,s-1}-\vr_{k-1,s}}
$$
from the space $\cS_{(k-1,s)}\oplus\cS_{(k-1,s-1)}$ to the first and second summand
respectively, where $\vr_{k,s}$ is the eigenvalue (\ref{rhoks}) of $\rR_0$.

{}From equation (\ref{NextLevel}) one gets
$$
\begin{array}{rcrl}
\displaystyle
P_{k-1,s}
&=&
\displaystyle
-\,\frac{\Pi^{k-1,s}_1}{\g_{k,s}-\g_{k-1,s}}\!&\!
\Big(
2(n\l+k-1)\rD(P_{k,s})-\rG_0(P_{k,s})
\Big)\\[14pt]
\displaystyle
P_{k-1,s-1}
&=&
\displaystyle
-\,\frac{\Pi^{k-1,s}_2}{\g_{k,s}-\g_{k-1,s-1}}\!&\!
\Big(
2(n\l+k-1)\rD(P_{k,s})-\rG_0(P_{k,s})
\Big).
\end{array}
$$
To rewrite the previous expression in terms of $\rD(P_{k,s})$ and $\rG_0(P_{k,s})$, one
resorts to the following formul{\ae}
$$
\begin{array}{rcl}
\rR_0\rD(P_{k,s})
&=&
\vr_{k,s}\rD(P_{k,s}) - 2\rG_0(P_{k,s})\\[8pt]
\rR_0\rG_0(P_{k,s})
&=&
2\rD(P_{k,s}) + 
\left(
\vr_{k,s}-2(n+2k-2)
\right)
\rG_0(P_{k,s})
\end{array}
$$
obtained with the help of the commutation relations of the generators (\ref{sl2}) and
(\ref{rh1}) of $\Sl(2,\bbR)\ltimes\rh_1$. A lengthy but straightforward calculation
gives
$$
P_{k-1,s}+P_{k-1,s-1}=\half\rD(P_{k,s})+\frac{(1-2\l)n}{2s(2s-2k-n+2)}\rG_0(P_{k,s}).
$$
Then, the definition (\ref{quantZeroShift}) of the quantization map yields the formula
(\ref{firstOrderIsomForm}).

(ii) In the harmonic case, $s=0$, the equation (\ref{NextLevel}) reduces to
\begin{equation}
\g_{k-1,0}P_{k-1,0}
=
\g_{k,0}P_{k-1,0}
+
2(n\l+k-1)\rD(P_{k,0})
\label{NextLevelHarm}
\end{equation}
since $\rG_0(P_{k,0})=0$. With the help of (\ref{gammaks}), one gets the formula
(\ref{firstOrderIsomFormHarm}).
\end{proof}

\begin{rmk}
{\rm
In the lower-order cases $k\leq2$, there exists an explicit formula for the
quantization map; it is given by the two formul{\ae} (\ref{Ansatz2},\ref{TheSolution})
and (\ref{QuantHalfDensity}) below.
}
\end{rmk}

With this preliminary result, we are ready to prove the announced theorem.

\begin{pro}
Given a differential linear operator $\cQ:\cS\to\cD_\l[[i\hbar]]$ of the form
$$
\cQ(P)=P+i\hbar(\a\rD(P)+\b\rG_0(P))+O(\hbar^2),
$$
the associative product $*:\cS\otimes\cS\to\cS$ defined by
$\cQ(P*Q)=\cQ(P)\circ\cQ(Q)$ is a star-product if and only if $\a=\half$ and $\b=0$.
\end{pro}
\begin{proof}
One can consider the inverse map $\cQ^{-1}:\cD_\l[[i\hbar]]\to\cS[[i\hbar]]$, which is given
by $\cQ^{-1}(P)=P-i\hbar(\a\rD(P)+\b\rG_0(P))+O(\hbar^2)$. Using the
well-known composition formula for differential operators
\begin{equation}
\begin{array}{rcl}
P\circ Q 
&= &
\displaystyle
\sum_{\ell=0}^\infty \frac{(i\hbar)^\ell}{\ell !}\, 
\partial_{\xi_{i_1}}\cdots\partial_{\xi_{i_\ell}}(P)\,
\partial_{i_1}\cdots\partial_{i_\ell}(Q),\\[14pt]
&=&
\displaystyle
PQ +i\hbar\partial_{\xi_j}(P)\partial_j(Q)+O(\hbar^2)
\end{array}
\label{composition}
\end{equation}
one obtains
$$
\begin{array}{rcl}
P*Q
&=&
\cQ^{-1}(\cQ(P)\circ\cQ(Q))\\[10pt]
&=&
PQ +i\hbar
\partial_{\xi_j}(P)\partial_j(Q)+
i\hbar\a
\Big(\rD(P)Q+P\rD(Q)-\rD(PQ)\Big)\\[10pt]
&&+
i\hbar\b
\Big(\rG_0(P)Q+P\rG_0(Q)-\rG_0(PQ)\Big) \\[10pt]
&&
+O(\hbar^2)\\[10pt]
&=&
\displaystyle
PQ +\frac{i\hbar}{\!2}\{P,Q\}
+
i\hbar\Big(\a-\half\Big)\Big(\rD(P)Q+P\rD(Q)-\rD(PQ)\Big)\\[10pt]
&&+
i\hbar\b
\Big(\rG_0(P)Q+P\rG_0(Q)-\rG_0(PQ)\Big) \\[10pt]
&&
+O(\hbar^2).
\end{array}
$$
Recall that the Hochschild boundary of a $1$-cochain~$A\in\End(\cS_\d)$ is given by
$(\rd{A})(P,Q)=A(P)Q+PA(Q)-A(PQ)$ and observe that the preceding expression is
therefore
\begin{equation}
P*Q
=
PQ +\frac{i\hbar}{\!2}\{P,Q\}
+
i\hbar\,\rd\Big(\Big(\a-\half\Big)D+\b\,\rG_0\Big)(P,Q)
+O(\hbar^2).
\end{equation}
One sees that $P*Q$ satisfies the definition (\ref{starDef}) of a
star-product if and only if $\a=\half$ and $\b=0$.
\end{proof}
The operation (\ref{starOperation}) is, actually, given by bi-differential
operators because the quantization map $\cQ_{\l;\hbar}$ given by (\ref{quantZeroShift})
and its inverse are differential operators at each order in~$\hbar$.
Indeed, we have $\widetilde\cQ_\l=\Id+\cN_\l$ as in Theorem \ref{ThmZeroShift}, so that
$(\cQ_{\l;\hbar})^{-1}$ is a differential operator as is $\cQ_{\l;\hbar}$.

Theorem \ref{starProduct} follows now from the preceding two propositions in the case
$k>2$ and from the explicit formula (\ref{QuantHalfDensity}) in the case $k\leq2$.

\section{Quantizing second-order polynomials}\label{quadratic}

This problem has first been solved in \cite{DO2}. It was proved
that  if $n=p+q\geq2$, there exists an isomorphism of $\so(p+1,q+1)$-modules
$
\widetilde{\cQ}^2_{\l,\m}
:
\cS^2_{\delta}\stackrel{\cong}{\longrightarrow}\cD^2_{\l,\m}
$
where $\delta=\m-\l$, provided
\begin{equation}
\delta\not\in\left\{\frac{2}{n},
\frac{n+2}{2n},
1,
\frac{n+1}{n},
\frac{n+2}{n}\right\}.
\label{Reson}
\end{equation}
This result is clearly consistent with the general Theorem \ref{IsomGen1}. Moreover,
the latter guarantees the uniqueness of such an isomorphisms under the further
condition that the principal symbol be preserved at each order.

\subsection{Explicit formul{\ae}}

In the non-resonant case, the explicit formula for the unique isomorphism has also been
computed in \cite{DO2}. One has
\begin{equation}
\widetilde{\cQ}^2_{\l,\m}
=
\Id+\gamma_1\rG_0
+\gamma_2\rD
+\gamma_3\Euler\rD
+\gamma_4\Delta_0
+\gamma_5\rD^2
\label{Ansatz2}
\end{equation}
where the numerical coefficients are given by
\begin{equation}
\matrix{
\displaystyle
\gamma_1=
\frac{n(\l+\m-1)}{2(n\delta-2)(n(\delta-1)-2)},
\hfill\cr\noalign{\medskip}
\displaystyle \gamma_2=
\frac{\l}{1-\delta},
\hfill\cr\noalign{\medskip}
\displaystyle \gamma_3=
\frac{1-\l-\m}{(\delta-1)(n(\delta-1)-2)},\hfill
\cr\noalign{\medskip}
\displaystyle \gamma_4=
\frac{n\l\Big(2+(4\l-1)n+(2\l^2-\l\m-\m^2+2\m-1)n^2\Big)}
{2(n(\delta-1)-1)(n(2\delta-1)-2)(n\delta-2)(n(\delta-1)-2)},
\hfill\cr\noalign{\medskip}
\displaystyle
\gamma_5=
\frac{n\l(n\l+1)}{2(n(\delta-1)-1)(n(\delta-1)-2)}.
\hfill\cr
}
\label{TheSolution}
\end{equation}

\goodbreak

In particular, the half-density quantization map (\ref{quantZeroShift}) is given by
\begin{equation}
\cQ^2_{\half;\hbar}
=
\Id
+\frac{i\hbar}{\!2}\rD
-\frac{\hbar^2}{\!8}\left(\frac{n}{(n+1)(n+2)}\Delta_0
+\frac{n}{(n+1)}\rD^2\right).
\label{QuantHalfDensity}
\end{equation}

\begin{rmk}
{\rm
At this stage, it is interesting to see how our conformally equivariant quantization
compares with the Weyl quantization on $T^*\bbR^n$. In our framework, the Weyl
quantization map, $\cQ_{\rm Weyl}$, retains the very elegant form
\begin{equation}
\begin{array}{rcl}
\displaystyle
\cQ_{\rm Weyl} 
&=&
\displaystyle\exp\Big(\frac{i\hbar}{\!2}\rD\Big)\\[10pt] &=&
\displaystyle
\Id
+\frac{i\hbar}{\!2}\rD
-\frac{\hbar^2}{\!8}\rD^2
+
O(\hbar^3)
\end{array}
\label{WeylMap}
\end{equation}
where the divergence operator $\rD$ is as in (\ref{rh1}).  (See, e.g., \cite{Fed}
p.~87.)
}
\end{rmk}

\goodbreak

\subsection{Study of the resonant modules}

For the sake of completeness, let us study in some more detail the particular modules
of differential operators corresponding to the resonances (\ref{Reson}). It has 
been shown~\cite{DO2} that, for each resonant value of $\delta$, there exist pairs
$(\l,\m)$ of weights such that the $\so(p+1,q+1)$-modules $\cS^2_\delta$ and
$\cD^2_{\l,\m}$ are isomorphic, namely
\begin{equation}
\setlength{\extrarowheight}{8pt}
\begin{array}{|c||c|c|c|c|c|}
\hline
\delta & \frac{2}{n} & \frac{n+2}{2n} & 1 & \frac{n+1}{n} &
\frac{n+2}{n}\\[8pt]
\hline
\hline
\lambda &
\frac{n-2}{2n} &
0,\frac{n-2}{2n} & 0 & 0, -\frac{1}{n} & -\frac{1}{n}\\[8pt]
\hline
\mu & \frac{n+2}{2n} & \frac{n+2}{2n},1 & 1 &
\frac{n+1}{n},1 & \frac{n+1}{n}\\[8pt]
\hline
\end{array}
\label{TheArray}
\end{equation}
However, in these cases, the isomorphism is not unique. For the particular values
$\d=2/n$, $1$, $(n+2)/n$, there is a unique choice of $(\l,\m)$ which, furthermore,
leads to symmetric quantized symbols; for example (see \cite{DO2}) the so-called
Yamabe operator (also know as the conformal Laplacian) shows up naturally in the
first resonant case in (\ref{TheArray}).

\subsection{Quantizing the geodesic flow}

Let us finally illustrate our quantization procedure with a specific and important
example, namely the quantization of the geodesic flow on a conformally flat
manifold~$(M,\rg)$.

Consider, on~$T^*M$, the quadratic Hamiltonian 
$$
H=\rg^{ij}\xi_i\xi_j.
$$ 
whose flow projects onto the geodesics of $(M,\rg)$.

Let us put $\l=\m=\half$ and apply, using (\ref{QuantHalfDensity}), the construction of
the quantum Hamiltonian (\ref{quantumOp}) spelled out in Section
\ref{quantumHamiltonian}. In doing so, we recover a result obtained in \cite{DO2},
namely
\begin{equation}
\hat{H}
=
-\hbar^2\left(
\Delta_\rg-\frac{n^2}{4(n-1)(n+2)}\,R_\rg
\right)
\label{QuantGeodFlow}
\end{equation}
where $R_\rg$ stands for the scalar curvature of $(M,\rg)$.
The operator (\ref{QuantGeodFlow}) is therefore a natural candidate for the quantum
Hamiltonian of the geodesic flow on a pseudo-Riemannian mani\-fold.




\end{document}